\documentclass[a4paper,12pt]{article}
\usepackage{latexsym,bm}
\usepackage{amsmath,amsfonts,amssymb}

\newcommand{\BS}{\mathfrak S}
\newcommand{\Z}{\mathbb Z} \newcommand{\HH}{\mathcal{H}}

\newcommand{\lam}{\lambda}
\newcommand{\ft}{\mathfrak t}

\newcommand{\eps}{\varepsilon}
\newcommand{\A}{\mathcal A}

\newcommand{\ksl}{\widehat{\mathfrak{sl}}}

\DeclareMathOperator{\GL}{GL} 
\DeclareMathOperator{\ch}{char} \DeclareMathOperator{\MM}{M}
\DeclareMathOperator{\cch}{ch} \DeclareMathOperator{\res}{res}
 \DeclareMathOperator{\Aut}{Aut}
\DeclareMathOperator{\rad}{rad} \DeclareMathOperator{\Irr}{Irr}
\DeclareMathOperator{\rres}{res}

\newtheorem{thm}{Theorem}\newtheorem{cor}{Corollary}
\newtheorem{lem}{Lemma}\newtheorem{dfn}{Definition}
\numberwithin{equation}{section} \numberwithin{prop}{section}
\numberwithin{thm}{section} \numberwithin{lem}{section}
\numberwithin{dfn}{section} \numberwithin{cor}{section}

\title{Mullineux involution and twisted affine Lie algebras
}
\author{Jun Hu$^{\star}$ \\[5pt]
$^{\star}$Department of Applied Mathematics\\
Beijing Institute of Technology\\
Beijing, 100081, P.R. China\null\\[1.5pt]
E-mail: junhu303@yahoo.com.cn}

\date{}

\begin{document}
\maketitle
\begin{abstract}

We use Naito-Sagaki's work [S. Naito \& D. Sagaki, J. Algebra 245
(2001) 395--412, J. Algebra 251 (2002) 461--474] on
Lakshmibai-Seshadri paths fixed by diagram automorphisms to study
the partitions fixed by Mullineux involution. We characterize the
set of Mullineux-fixed partitions in terms of crystal graphs of
basic representations of twisted affine Lie algebras of type
$A_{2\ell}^{(2)}$ and of type $D_{\ell+1}^{(2)}$. We set up
bijections between the set of symmetric partitions and the set of
partitions into distinct parts. We propose a notion of double
restricted strict partitions. Bijections between the set of
restricted strict partitions (resp., the set of double restricted
strict partitions) and the set of Mullineux-fixed partitions in
the odd case (resp., in the even case) are obtained.
\end{abstract}
\bigskip\bigskip\bigskip

\section{Introduction}

Let $n, e\in\mathbb{N}$. Let $k$ be a field and $0\neq q\in k$.
Suppose that either $e>1$ and $q$ is a primitive $e$th root of
unity; or $q=1$ and $\ch k=e$.\footnote{In the latter case, $e$ is
necessarily to be a prime number.} Let $\HH_k(\BS_n)$ be the
Iwahori-Hecke algebra associated to the symmetric group $\BS_n$
with parameter $q$ and defined over $k$. The Mullineux involution
$\MM$ is a bijection defined on the set of all $e$-regular
partitions of $n$, which arises naturally when one twists
irreducible modules (labelled by $e$-regular partitions) over
$\HH_k(\BS_n)$ by a $k$-algebra automorphism $\#$ (see Section 2
for definition of $\#$). If $q=1$ and $e$ is an odd prime number,
the involution $\MM$ determines which simple module splits and
which remains simple when restricting to the alternating subgroup
$A_n$. In that case, the set of partitions which are fixed by the
involution $\MM$ parameterizes the irreducible modules of $k\BS_n$
which split on restriction to $A_n$. In \cite{Kl}, Kleshchev gave
a remarkable algorithm for computing the involution $\MM$. A
crystal bases approach to Kleshchev's algorithm of the involution
$\MM$ was given in \cite[(7.1)]{LLT}.\smallskip

The purpose of this paper is to study the partitions fixed by
Mullineux involution for arbitrary $e$. We find that the set of
Mullineux-fixed partitions is related to the twisted affine Lie
algebras of type $A_{2\ell}^{(2)}$ and of type $D_{\ell+1}^{(2)}$,
which reveals new connection between the theory of affine Lie
algebra and the theory of modular representations. Our main tool
are Naito-Sagaki's work (\cite{NS2}, \cite{NS1}) on
Lakshmibai-Seshadri paths fixed by diagram automorphisms, which
was also used in \cite{Hu4} and \cite{Hu5} to derive explicit
formulas for the number of modular irreducible representations of
the cyclotomic Hecke algebras of type $G(r,p,n)$, see \cite{Hu1},
\cite{Hu2} and \cite{Hu3} for related work. We characterize the
set of Mullineux-fixed partitions in terms of crystal graph of
basic representations of twisted affine Lie algebras of type
$A_{2\ell}^{(2)}$ and of type $D_{\ell+1}^{(2)}$ (Theorem
\ref{thm37}). We set up bijections (Theorem \ref{thm315} and
Theorem \ref{thm317}) between the set of Mullineux-fixed
partitions in the odd case (resp., the set of symmetric
partitions) and the set of restricted strict partitions (resp.,
the set of partitions into distinct parts). As an application, we
obtain new identities on the cardinality of the set of
Mullineux-fixed partitions in terms of the principal specialized
characters of the basic representations of these twisted affine
Lie algebras (Theorem \ref{thm313} and Theorem \ref{thm320}).
Furthermore, we propose a notion of double restricted strict
partitions (Definition \ref{dfn321}), which is a direct explicit
characterization of Kang's reduced proper Young wall of type
$D_{\ell+1}^{(2)}$ (\cite{K}). We obtain a bijection (Theorem
\ref{thm324}) between the set of Mullineux-fixed partitions in the
even case and the set of double restricted strict partitions. Our
main results shed some new insight on the modular representations
of the alternating group and of Hecke-Clifford superalgebras as
well as of the spin symmetric group (see Remark 3.25 and Remark
3.18), which clearly deserves further study.
\bigskip\bigskip\bigskip

\section{Preliminaries}

In this section, we shall first review some basic facts about the
representation of the Iwahori-Hecke algebras associated to
symmetric groups. Then we shall introduce the notion of Mullineux
involution, Kleshchev's $e$-good lattice as well as Kleshchev's
algorithm of Mullineux involution.
\smallskip

Let $\BS_n$ be the symmetric group on $\{1,2,\cdots,n\}$, acting
from the right. Let ${\A}=\Z[v,v^{-1}]$, where $v$ is an
indeterminate. The Iwahori-Hecke algebra $\HH_{\A}(\BS_n)$
associated to $\BS_n$ is the associative unital ${\A}$-algebra
with generators $T_1,\cdots,$ $T_{n-1}$ subject to the following
relations
$$\begin{aligned}
&(T_i-v)(T_i+1)=0,\quad\text{for $1\leq i\leq n-1$,}\\
&T_iT_{i+1}T_i=T_{i+1}T_{i}T_{i+1},\quad\text{for $1\leq i\leq n-2$,}\\
&T_iT_j=T_jT_i,\quad\text{for $1\leq i<j-1\leq n-2$.}\end{aligned}
$$
For each integer $i$ with $1\leq i\leq n-1$, we define
$s_i=(i,i+1)$. Then $S:=\{s_1,s_2,\cdots,s_{n-1}\}$ is the set of
all the simple reflections in $\BS_n$. A word $w=s_{i_1}\cdots
s_{i_k}$ for $w\in\BS_{n}$ is a reduced expression if $k$ is
minimal; in this case we say that $w$ has length $k$ and we write
$\ell(w)=k$. Given a reduced expression $s_{i_1}\cdots s_{i_k}$
for $w\in\BS_n$, we write $T_w=T_{i_1}\cdots T_{i_k}$. The braid
relations for generators $T_1,\cdots,T_{n-1}$ ensure that $T_w$ is
independent of the choice of reduced expression. It is well-known
that $\HH_{{\A}}(\BS_n)$ is a free ${\A}$-module with basis
$\{T_w|w\in\BS_n\}$. For any field $k$ which is an ${\A}$-algebra,
we define $\HH_k(\BS_n):=\HH_{{\A}}(\BS_n)\otimes_{{\A}}k$. Then
$\HH_k(\BS_n)$ can be naturally identified with the $k$-algebra
defined by the same generators and relations as $\HH_{\A}(\BS_n)$
above. Specializing $v$ to $1\in k$, one recovers the group
algebra $k\BS_n$ of $\BS_n$ over $k$.\smallskip

We recall some combinatorics. A partition of $n$ is a
non-increasing sequence of positive integers
$\lam=(\lam_1,\cdots,\lam_r)$ such that $\sum_{i=1}^{r}\lam_i=n$.
For any partition $\lam=(\lam_1,\lam_2,\cdots)$, the conjugate of
$\lam$ is defined to be a partition
$\lam^t=(\lam_1^t,\lam_2^t,\cdots)$, where
$\lam_j^t:=\#\{i|\lam_i\geq j\}$ for $j=1,2,\cdots$. We define
$\ell(\lam):=\max\{i|\lam_i\neq 0\}$. For any partition $\lam$ of
$n$,  we denote by $\ft^{\lam}$ (resp., $\ft_{\lam}$) the standard
$\lam$-tableau in which the numbers $1,2,\cdots,n$ appear in order
along successive rows (resp., columns). The row stabilizer of
$\ft^{\lam}$, denoted by $\BS_{\lam}$, is the standard Young
subgroup of $\BS_n$ corresponding to $\lam$. Let
$$
x_{\lam}=\sum_{w\in\BS_{\lam}}T_w,\quad
y_{\lam}=\sum_{w\in\BS_{\lam}}(-v)^{-\ell(w)}T_w. $$ Let
$w_{\lam}\in\BS_n$ be such that $\ft^{\lam}w_{\lam}=\ft_{\lam}$.
Following \cite[Section 4]{DJ1}, we define
$z_{\lam}=x_{\lam}T_{w_{\lam}}y_{\lam^t}$.

\begin{dfn}\label{df21} The right ideal $z_{\lam}\HH$ is called the right Specht module
of $\HH=\HH_{{\A}}(\BS_n)$ corresponding to $\lam$. We denote it
by $S^{\lam}$.
\end{dfn}

For any field $k$ which is an ${\A}$-algebra, let
$S_k^{\lam}:=S^{\lam}\otimes_{\A} k$. There is a natural bilinear
form ${\langle,\rangle}$ on each $S^{\lam}$ (and hence on each
$S_k^{\lam}$). Let $D_k^{\lam}:=S_k^{\lam}/\rad\langle,\rangle$.
Let ``$\trianglelefteq$'' be the dominance order on the set of all
partitions as defined in \cite[(3.1)]{Mu}.

\addtocounter{lem}{1}
\begin{lem} {\rm (\cite{DJ1})}\label{lm1} With the above notations, we have
\smallskip

1) the set of all the nonzero $D_k^{\lam}$ (where $\lam$ runs over
partitions of $n$) forms a complete set of pairwise non-isomorphic
simple $\HH_{k}(\BS_{n})$-modules. Moreover, if $\HH_{k}(\BS_{n})$
is semisimple, then $D_k^{\lam}=S_k^{\lam}\neq 0$ for every
partition $\lam$ of $n$;
\smallskip

2) if $D_k^{\mu}\neq 0$ is a composition factor of $S_k^{\lambda}$
then $\lam\trianglelefteq\mu$, and every composition factor of
$S_k^{\lam}$ is isomorphic to some $D_k^{\mu}$ with
$\lam\trianglelefteq\mu$. If $D_k^{\lam}\neq 0$ then the
composition multiplicity of $D_k^{\lam}$ in $S_k^{\lam}$ is $1$.
\end{lem}

Henceforth, let $k$ be a fixed field which is an ${\A}$-algebra.
We assume that $v$ is specialized to $q\in k$ such that
$1+q+q^2+\cdots+q^{a-1}=0$ for some positive integer $a$. We
define
$$ e=\min\bigl\{1<a<\infty\bigm|\text{$1+q+q^2+\cdots+q^{a-1}=0$
in $k$}\bigr\}.
$$
Clearly, $e=\ch k$ if $q=1$; and otherwise $e$ is the
multiplicative order of $q$. For simplicity, we shall write
$\HH_k$ instead of $H_k(\BS_n)$.

A partition $\lam$ is called $e$-regular if it contains at most
$e-1$ repeating parts, i.e., $\lam=(1^{m_1}2^{m_2}\cdots
j^{m_j}\cdots)$ with $0\leq m_i<e$ for every $i$. By \cite{DJ1},
for any partition $\lam$ of $n$, $D_k^{\lam}\neq 0$ if and only if
$\lam$ is $e$-regular. Let $\mathcal{K}_n$ be the set of all the
$e$-regular partitions of $n$. Let $\#$ (see \cite{DJ1},
\cite[(2.3)]{Mu}) be the $k$-algebra automorphism of $\HH_k$ which
is defined on generators by $T_{i}^{\#}=-vT_{i}^{-1}$ for each
$1\leq i<n$. For each $\HH_k(\BS_n)$-module $V$, we denote by
$V^{\#}$ the $\HH_k(\BS_n)$-module obtained by twisting $V$ by
$\#$. That is, $V^{\#}=V$ as $k$-linear space, and $v\cdot
h:=vh^{\#}$ for any $v\in V$ and $h\in \HH_k(\BS_n)$. Let $\ast$
be the algebra anti-automorphism on $\HH_k$ which is defined on
generators by $T_i^{\ast}=T_{i}$ for any $1\leq i<n$.

\addtocounter{dfn}{1}
\begin{dfn} {\rm (\cite{Mul}, \cite{Br})} Let $\MM$ be the unique involution defined on the set
$\mathcal{K}_n$ such that $\bigl(D_k^{\lam}\bigr)^{\#}\cong
D_k^{\MM(\lam)}$ for any $\lam\in\mathcal{K}_n$. We call the map
$\MM$ the Mullineux involution, and $\lam$ a Mullineux-fixed
partition if $\MM(\lam)=\lam$.
\end{dfn}

An algorithm which compute the involution $\MM$ was first proposed
by Mullineux in 1979, when he constructed an involution on the set
of $e$-regular partitions and conjectured its coincidence with the
above $\MM$. Mullineux worked in the setup that $q=1$ and $e$
being a prime number, though his combinatorial algorithm does not
really depend on $e$ being prime. In \cite{Kl}, Kleshchev gave a
quite different remarkable algorithm of the involution $\MM$ based
on his work of branching rules for the modular representations of
symmetric groups. In \cite{FK}, Ford and Kleshchev proved that
Kleshchev's algorithm is equivalent to Mullineux's original
algorithm and thus proved Mullineux's conjecture. The validity of
Kleshchev's algorithm of $\MM$ for arbitrary $e$ is proved in
\cite{Br}.\smallskip

Note that the Mullineux involution $\MM$ depends only on $e$. {\it
Henceforth, we refer to the case when $e$ is odd as the odd case;
and to the case when $e$ is even as the even case.} By
\cite[(3.5)]{DJ2} and \cite[(5.2),(5.3)]{Mu},
$\bigl(S^{\lam}\bigr)^{\#}\cong\bigl(S^{\lam^t}\bigr)^{\ast}$. If
$\HH_{k}(\BS_{n})$ is semisimple, then
$\bigl(S_k^{\lam^t}\bigr)^{\ast}\cong S_k^{\lam^t}$, hence in that
case the involution $\MM$ degenerates to the map
$\lam\mapsto\lam^t$. In this paper, we do not need Mullineux's
original combinatorial algorithm (\cite{Mul}) for defining $\MM$,
but we do need Kleshchev's algorithm (\cite{Kl}) of the involution
$\MM$. To this end, we have to recall the notion of Kleshchev's
$e$-good lattice.\smallskip

Let $\lam$ be a partition of $n$. The Young diagram of $\lam$ is
the set
$$ [\lam]=\bigl\{(a,b)\bigm|\text{$1\leq b\leq\lam_{a}$}\bigr\}. $$
The elements of $[\lam]$ are nodes of $\lam$. Given any two nodes
$\gamma=(a,b), \gamma'=(a',b')$ of $\lam$, say that $\gamma$ is
{\it below} $\gamma'$, or $\gamma'$ is {\it above} $\gamma$, if
$a>a'$. The {\it residue} of $\gamma=(a,b)$ is defined to be
$\rres(\gamma):=b-a+e\Z\in\Z/e\Z$, and we say that $\gamma$ is a
$\rres(\gamma)$-node. Note that we can identify the set
$\{0,1,2,\cdots,e-1\}$ with $\Z/e\Z$ via $i\mapsto \overline{i}$
for each $0\leq i\leq e-1$. Therefore, we can also think that the
$\res(?)$ function takes values in $\{0,1,2,\cdots,e-1\}$.

A removable node is a node of the boundary of the Young diagram
$[\lam]$ which can be removed, while an addable node is a concave
corner on the rim of $[\lam]$ where a node can be added. If $\mu$
is a partition of $n+1$ with $[\mu]=[\lam]\cup
\bigl\{\gamma\bigr\}$ for some removable node $\gamma$ of $\mu$,
we write $\lam\rightarrow\mu$. If in addition $\res(\gamma)=x$, we
also write that $\lam\overset{x}{\rightarrow}\mu$. For example,
suppose $n=42$ and $e=3$. The nodes of $\lam=(9^2,8,7,5,3,1)$ have
the following residues $$ \lam=\left(\begin{matrix} \overline{0} &
\overline{1}& \overline{2}& \overline{0}& \overline{1} &
\overline{2} & \overline{0}
& \overline{1} & \overline{2}   \\
\overline{2} & \overline{0} & \overline{1} & \overline{2} &
\overline{0} & \overline{1} & \overline{2} &\overline{0} &  \overline{1}  \\
\overline{1}& \overline{2}& \overline{0}& \overline{1}&
\overline{2}& \overline{0} & \overline{1}& \overline{2}&  \\
\overline{0}& \overline{1}& \overline{2}& \overline{0}&
\overline{1}& \overline{2} & \overline{0}& & \\
\overline{2}& \overline{0} & \overline{1}& \overline{2}
&\overline{0} &&&& \\
\overline{1}& \overline{2} & \overline{0} & & & & & & \\
\overline{0}&  & & & & & & &
\end{matrix}
\right) .
$$
It has six removable nodes. Fix a residue $x$ and consider the
sequence of removable and addable $x$-nodes obtained by reading
the boundary of $\lam$ from the top down. In the above example, if
we consider residue $x=\overline{0}$, then we get a sequence
AARRRR, where each ``A'' corresponds to an addable
$\overline{0}$-node and each ``R'' corresponds to a removable
$\overline{0}$-node. Given such a sequence of letters A,R, we
remove all occurrences of the string ``AR'' and keep on doing this
until no such string ``AR'' is left. The ``R''s that still remain
are the {\it normal} $\overline{0}$-nodes of $\lam$ and the
rightmost of these is the {\it good} $\overline{0}$-node. In the
above example, the two removable $\overline{0}$-nodes in the last
two rows survive after we delete all the string ``AR''. Therefore,
the removable $\overline{0}$-node in the last row is the good
$\overline{0}$-node. If $\gamma$ is a good $x$-node of $\mu$ and
$\lam$ is the partition such that $[\mu]=[\lam]\cup\gamma$, we
write $\lam\overset{x}{\twoheadrightarrow}\mu$. The {\it
Kleshchev's $e$-good lattice} is, by definition, the infinite
graph whose vertices are the $e$-regular partitions and whose
arrows are given by $
\text{$\lam\overset{x}{\twoheadrightarrow}\mu$\,\,\,$\Longleftrightarrow$\,
$\lam$ is obtained from $\mu$ by removing a good $x$-node}$. It is
well-known that, for each $e$-regular partition $\lam$, there is a
path (not necessary unique) from the empty partition $\emptyset$
to $\lam$ in Kleshchev's $e$-good lattice.\medskip

Kleshchev's $e$-good lattice in fact provides a combinatorial
realization of the crystal graph of the basic representation of
the affine Lie algebra of type $A_{e-1}^{(1)}$ (which we denote by
$\ksl_{e}$). To be more precise, let
$\{\alpha_0,\alpha_1,\cdots,\alpha_{e-1}\}$ be the set of simple
roots of $\ksl_{e}$, let
$\bigl\{\alpha_0^{\vee},\alpha_1^{\vee},\cdots,\alpha_{e-1}^{\vee}\bigr\}$
be set of simple coroots, let $$
\begin{pmatrix}
2& -1& 0& \cdots & 0& -1\\
-1& 2& -1& \cdots & 0& 0\\
0& -1& 2& \cdots & 0& 0\\
\vdots & \vdots & \vdots & &\vdots &\vdots \\
0& 0& 0& \cdots & 2& -1\\
-1& 0& 0& \cdots & -1& 2
\end{pmatrix}_{e\times e}\quad\text{if $e\geq 3$;}
$$ or $$\begin{pmatrix}
2& -2\\
-2& 2
\end{pmatrix}_{2\times 2}\quad\text{if $e=2.$}
$$
be the corresponding affine Cartan matrix. Let $d$ be the scaling
element. Then the set
$\bigl\{\alpha_0^{\vee},\alpha_1^{\vee},\cdots,\alpha_{e-1}^{\vee},d\bigr\}$
forms a basis of the Cartan subalgebra of $\ksl_{e}$, let $\bigl\{
\Lambda_0,\Lambda_1,\cdots,\Lambda_{e-1},\delta\bigr\}$ be the
corresponding dual basis, where $\delta$ denotes the null root.
The integrable highest weight module of highest weight
$\Lambda_0$, denoted by $L(\Lambda_0)$, is called the basic
representation of $\ksl_{e}$. It is a remarkable fact (\cite{MM},
\cite[(2.11)]{AM}) that the crystal graph of $L(\Lambda_0)$ is
exactly the same as the Kleshchev's $e$-good lattice if one use
the embedding $L(\Lambda_0)\subset\mathcal{F}(\Lambda_0)$, where
$\mathcal{F}(\Lambda_0)$ is the Fock space as defined in
\cite[\S4.2]{LLT}. In particular, an explicit formula for the
number of irreducible $\HH_k(\BS_n)$-modules, i.e.,
$\#\mathcal{K}_n$, is known (see \cite{AM}), which was expressed
in terms of principal specialized character of the basic
representation $L(\Lambda_0)$.\smallskip

Now we can state Kleshchev's algorithm of the Mullineux involution
$\MM$. Here we follow Lascoux-Lerclerc-Thibon's reformulation in
\cite[(7.1)]{LLT}.

\addtocounter{lem}{1}
\begin{lem} \label{lm24} {\rm (\cite{Kl})} Let $\lam\in\mathcal{K}_n$ be an $e$-regular partition of $n$,
and let $$ \emptyset\overset{r_1}{\twoheadrightarrow}\cdot
\overset{r_2}{\twoheadrightarrow}\cdot \cdots\cdots
\overset{r_n}{\twoheadrightarrow}\lam $$ be a path from
$\emptyset$ to $\lam$ in Kleshchev's $e$-good lattice. Then, the
sequence $$
\underline{\emptyset}\overset{e-r_1}{\twoheadrightarrow}\cdot
\overset{e-r_2}{\twoheadrightarrow}\cdot \cdots\cdots
\overset{e-r_n}{\twoheadrightarrow}\cdot
$$
also defines a path in Kleshchev's $e$-good lattice, and it
connects $\emptyset$ to $\MM(\lam)$.
\end{lem}

Note that the Mullineux involution $\MM$ gives rise to an
equivalence relation on ${\mathcal K}_n$. That is, $\lam\sim\mu$
if and only if $\lam=\MM(\mu)$ for any $\lam,\mu\in\mathcal{K}_n$.
Let $A_n$ be the alternating group, which is a normal subgroup in
$\BS_n$ of index $2$. In the special case where $q=1$ and $e$ is
an odd prime number, the involution $\MM$ is closely related to
the modular representation of the alternating group $A_n$, as can
be seen from the following lemma.

\begin{lem} {\rm (\cite[(2.1)]{F})} Suppose that $q=1$ and $e$ is an odd prime number.
In particular, $\ch k=e$. Assume that $A_n$ is split over $k$.
Then

(1) for any $\lam\in\mathcal{K}_n$ with $\MM(\lam)\neq\lam$,
$D^{\lam}\downarrow_{A_n}$ remains irreducible;

(2) for any $\lam\in\mathcal{K}_n$ with $\MM(\lam)=\lam$,
$D^{\lam}\downarrow_{A_n}$ is a direct sum of two irreducible,
non-equivalent, representations of $kA_n$, say $D_+^{\lam}$ and
$D_{-}^{\lam}$;

(3) the set $$ \Bigl\{D_{+}^{\lam},
D_{-}^{\lam}\Bigm|\lam\in\mathcal{K}_n/{\sim},
\MM(\lam)=\lam\Bigl\} \bigsqcup
\Bigl\{D^{\lam}\downarrow_{A_n}\Bigm|\lam\in\mathcal{K}_n/{\sim},
\MM(\lam)\neq\lam\Bigr\}
$$
forms a complete set of pairwise non-isomorphic irreducible
$kA_n$-modules.
\end{lem}

As a consequence, we get that $$\begin{aligned}
&\quad\,\#\Irr\bigl(kA_n\bigr)\\
&=\frac{1}{2}\Bigl(\#\mathcal{K}_n-\#\bigl\{\lam\in\mathcal{K}_n\bigm|\MM(\lam)=\lam\bigr\}
\Bigr)+2\#\bigl\{\lam\in\mathcal{K}_n\bigm|\MM(\lam)=\lam\bigr\}\\
&=\frac{1}{2}\Bigl(\#\mathcal{K}_n+3\#\bigl\{\lam\in\mathcal{K}_n\bigm|\MM(\lam)=\lam\bigr\}
\Bigr).\end{aligned}
$$
\bigskip\bigskip\bigskip

\section{The orbit Lie algebras}

In this section, we shall first determine the orbit Lie algebras
corresponding to the Dynkin diagram automorphisms arising from the
Mullineux involution. Then we shall use Naito-Sagaki's work
(\cite{NS2}, \cite{NS1}) to study the set of Mullineux-fixed
partitions in terms of crystal graphs of basic representations of
the orbit Lie algebras, which are some twisted affine Lie algebras
of type $A_{2\ell}^{(2)}$ or of type $D_{\ell+1}^{(2)}$. The main
results are given in Theorem \ref{thm37}, Theorem \ref{thm313},
Theorem \ref{thm315}, Theorem \ref{thm317}, Theorem \ref{thm320}
and Theorem \ref{thm324}.
\medskip

Let $\mathfrak{g}$ be the Kac-Moody algebra over $\mathbb{C}$
associated to a symmetrizable generalized Cartan matrix
$(a_{i,j})_{i,j\in I}$ of finite size, where
$I=\{0,1,\cdots,e-1\}$. Let $\mathfrak{h}$ be its Cartan
subalgebra, and $W$ be its Weyl group. Let
$\{\alpha_i^{\vee}\}_{0\leq i\leq e-1}$ be the set of simple
coroots in $\mathfrak{h}$. Let
$\mathcal{X}:=\bigl\{\Lambda\in\mathfrak{h}^{\ast}\bigm|
\Lambda(\alpha_i^{\vee})\in\Z,\,\forall\,0\leq i<e\bigr\}$ be the
weight lattice. Let
$\mathcal{X}^{+}:=\bigl\{\Lambda\in\mathcal{X}\bigm|\Lambda(\alpha_i^{\vee})\geq
0,\,\forall\,0\leq i<e\bigr\}$ be the lattice of integral dominant
weights. Let
$\mathcal{X}_{\mathbb{R}}:=\mathcal{X}\otimes_{\Z}\mathbb{R}$,
where $\mathbb{R}$ is the field of real numbers. Assume that
$\Lambda\in\mathcal{X}^{+}$. P. Littelmann introduced (\cite{Li1},
\cite{Li2}) the notion of Lakshmibai-Seshadri paths (L-S paths for
short) of class $\Lambda$, which are piecewise linear, continuous
maps $\pi:[0,1]\rightarrow\mathcal{X}_{\mathbb{R}}$ parameterized
by pairs $(\underline{\nu},\underline{a})$ of a sequence
$\underline{\nu}: \nu_1>\nu_2>\cdots>\nu_s$ of elements of
$W\Lambda$, where $>$ is the ``relative Bruhat order" on
$W\Lambda$, and a sequence $\underline{a}: 0=a_0<a_1<\cdots<a_s=1$
of rational numbers with a certain condition, called the chain
condition. The set $\mathbb{B}(\Lambda)$ of all L-S paths of class
$\Lambda$ is called the path model for the integrable highest
weight module $L(\Lambda)$ of highest weight $\Lambda$ over
$\mathfrak{g}$. It is a remarkable fact that $\mathbb{B}(\Lambda)$
has a canonical crystal structure isomorphic to the crystal (in
the sense of \cite{Kas}) associated to the integrable highest
weight module of highest weight $\Lambda$ over the quantum algebra
$U'_v(\mathfrak{g})$ .\smallskip

Now let $\mathfrak{g}$ be the affine Kac-Moody algebra of type
$A_{e-1}^{(1)}$. Let $\omega:\,I\rightarrow I$ be an involution
defined by $\omega(0)=0$ and $\omega(i)=e-i$ for any $0\neq {i}\in
I$.

\begin{lem} $\omega$ is a Dynkin diagram automorphism in the sense of
\cite[\S1.2]{NS1}. That is $a_{\omega(i),\omega(j)}=a_{i,j}$,
$\forall\,i,j\in I$.
\end{lem}

\noindent {Proof:} \,This follows from direct verification.
\medskip

By \cite{FSS}, $\omega$ induces a Lie algebra automorphism (which
are called diagram outer automorphism)
$\omega\in\Aut(\mathfrak{g})$ of order $2$ and a linear
automorphism $\omega^{\ast}\in\GL(\mathfrak{h}^{\ast})$ of order
$2$. Following \cite{FRS} and \cite[\S1.3]{NS1} (where they work
with an arbitrary Kac-Moody algebra $\mathfrak{g}$ and a Dynkin
diagram automorphism $\omega$), we set
$c_{i,j}:=\sum\limits_{t=0}^{N_j-1}a_{i,\omega^t(j)}$, where
$N_j:=\#\bigl\{\omega^t(i)\bigm|t\geq 0\bigr\}$, $i,j\in I$. We
choose a complete set $\widehat{I}$ of representatives of the
$\omega$-orbits in $I$, and set
$\check{I}:=\bigl\{i\in\widehat{I}\bigm|c_{i,i}>0\bigr\}$. We put
$\hat{a}_{i,j}:=2c_{i,j}/c_{j}$ for $i,j\in\widehat{I}$, where
$c_i:=c_{ii}$ if $i\in\check{I}$, and $c_i:=2$ otherwise. Then
$(\hat{a}_{i,j})_{i,j\in\widehat{I}}$ is a symmetrizable
Borcherds-Cartan matrix (\cite{Bo}), and (if
$\check{I}\neq\emptyset$) its submatrix
$(\hat{a}_{i,j})_{i,j\in\check{I}}$ is a generalized Cartan
matrix. Let $\widehat{\mathfrak{g}}$ be the generalized Kac-Moody
algebra over $\mathbb C$ associated to
$(\hat{a}_{i,j})_{i,j\in\widehat{I}}$, with Cartan subalgebra
$\widehat{\mathfrak{h}}$, Chevalley generators
$\{\hat{x}_i,\hat{y}_i\}_{i\in\widehat{I}}$. The orbit Lie algebra
$\check{\mathfrak{g}}$ is defined to be the subalgebra of
$\widehat{\mathfrak{g}}$ generated by $\widehat{\mathfrak{h}}$ and
$\hat{x}_i,\hat{y}_i$ for $i\in\check{I}$, which is a usual
Kac-Moody algebra.

\begin{lem} With the above assumptions and notations, we have that
in our special case, $\check{\mathfrak{g}}$ is isomorphic to the
twisted affine Lie algebra of type $A_{2\ell}^{(2)}$ if
$e=2\ell+1$; and $\check{\mathfrak{g}}$ is isomorphic the twisted
affine Lie algebra of type $D_{\ell+1}^{(2)}$ if $e=2\ell$.
\end{lem}

\noindent {Proof:} \, We divide the proof into two cases:
\smallskip

\noindent {\it Case 1.}\,\,$e=2\ell+1$. The involution $\omega$ is
given by $$ \omega:\left\{\begin{aligned}{0}&\mapsto
{0}\\
{1}&\mapsto
{2\ell}\\
& \vdots\\
{\ell-1}&\mapsto {\ell+2}\\
{\ell}&\mapsto {\ell+1}
\end{aligned}\right. ,\quad\,\,
\left\{\begin{aligned} {\ell+1}&\mapsto {\ell}\\
& \vdots\\
{2\ell-1}&\mapsto {2}\\
{2\ell}&\mapsto {1}
\end{aligned}\right. .
$$
It is easy to check that $c_{i,i}=2$ for any $0\leq i<\ell$ and
$c_{\ell,\ell}=1$. We shall take
$\widehat{I}=\{{0},{1},\cdots,{l}\}$. By direct verification, we
get that $\check{I}=\widehat{I}$ and $$
(\hat{a}_{i,j})_{i,j\in\widehat{I}}=\begin{pmatrix}
2& -2& 0& \cdots & 0& 0& 0\\
-1& 2& -1& \cdots & 0& 0& 0\\
0& -1& 2& \cdots & 0& 0& 0\\
\vdots & \vdots & \vdots & & \vdots &\vdots &\vdots \\
0& 0& 0& \cdots & 2& -1& 0\\
0& 0& 0& \cdots & -1& 2& -2\\
0& 0& 0& \cdots & 0& -1& 2
\end{pmatrix}_{(\ell+1)\times (\ell+1)}\quad\text{if $\ell\geq 2$;}
$$ or $$\begin{pmatrix}
2& -4\\
-1& 2
\end{pmatrix}_{2\times 2}\quad\text{if $\ell=1.$}$$
Clearly this is an affine Cartan matrix of type $A_{2\ell}^{(2)}$,
hence in this case $\check{\mathfrak{g}}$ is isomorphic to the
twisted affine Lie algebra of type $A_{2\ell}^{(2)}$.\medskip

\noindent {\it Case 2.}\,\,$e=2\ell$. The involution $\omega$ is
given by $$ \omega:\left\{\begin{aligned}{0}&\mapsto
{0}\\
{1}&\mapsto
{2\ell-1}\\
& \vdots\\
{\ell-1}&\mapsto {\ell+1}\\
{\ell}&\mapsto{\ell}
\end{aligned}\right.,\quad\,\,\left\{\begin{aligned}{\ell+1}&\mapsto {\ell-1}\\
& \vdots\\
{2\ell-2}&\mapsto
{2}\\
{2\ell-1}&\mapsto {1}
\end{aligned}\right. .
$$
It is easy to check that $c_{i,i}=2$ for any $0\leq i\leq\ell$. We
shall take $\widehat{I}=\{{0},{1},\cdots,{l}\}$. By direct
verification, we get that $\check{I}=\widehat{I}$ and $$
(\hat{a}_{i,j})_{i,j\in\widehat{I}}=\begin{pmatrix}
2& -2& 0& \cdots & 0& 0& 0\\
-1& 2& -1& \cdots & 0& 0& 0\\
0& -1& 2& \cdots & 0& 0& 0\\
\vdots & \vdots & \vdots & & \vdots &\vdots &\vdots \\
0& 0& 0& \cdots & 2& -1& 0\\
0& 0& 0& \cdots & -1& 2& -1\\
0& 0& 0& \cdots & 0& -2& 2
\end{pmatrix}_{(\ell+1)\times (\ell+1)}\quad\text{if $\ell\geq 2$;}
$$ or $$\begin{pmatrix}
2& -2\\
-2& 2
\end{pmatrix}_{2\times 2}\quad\text{if $\ell=1.$}$$
Clearly this is an affine Cartan matrix of type
$D_{\ell+1}^{(2)}$, hence in this case $\check{\mathfrak{g}}$ is
isomorphic to the twisted affine Lie algebra of type
$D_{\ell+1}^{(2)}$.\medskip

We define
$\bigl(\mathfrak{h}^{\ast}\bigr)^{\circ}:=\bigl\{\Lambda\in\mathfrak{h}^{\ast}
\bigm|\omega^{\ast}(\Lambda)=\Lambda\bigr\}$.
$\widetilde{W}:=\bigl\{w\in
W\bigm|\omega^{\ast}w=w\omega^{\ast}\bigr\}$. We indicate by
$\check{}$ the objects for the obit Lie algebra
$\check{\mathfrak{g}}$. For example, $\check{\mathfrak{h}}$
denotes the Cartan subalgebra of $\check{\mathfrak{g}}$,
$\check{W}$ the Weyl group of $\check{\mathfrak{g}}$,
$\{\check{\Lambda}_i\}_{0\leq i\leq \ell}$ the set of fundamental
dominant weights in $\check{\mathfrak{h}}^{\ast}$. There exists a
linear automorphism
$P_{\omega}^{\ast}:\,\check{\mathfrak{h}}^{\ast}\rightarrow
\bigl(\mathfrak{h}^{\ast}\bigr)^{\circ}$ and a group isomorphism
$\Theta:\,\check{W}\rightarrow\widetilde{W}$ such that
$\Theta(\check{w})=P_{\omega}^{\ast}\check{w}\bigl(P_{\omega}^{\ast}\bigr)^{-1}$
for each $w\in\check{W}$. By \cite[\S6.5]{FSS}, for each $0\leq
i\leq \ell$, $$
P_{\omega}^{\ast}(\check{\Lambda}_i)=\sum_{t=0}^{N_i-1}\Lambda_{\omega^t(i)}+C\delta,
$$
where $N_i$ denotes the number of elements in the $\omega$-orbit
of $i$, $C\in\mathbb{Q}$ is some constant depending on $\omega$,
$\delta$ denotes the null root of $\mathfrak{g}$. It follows that
$P_{\omega}^{\ast}(\check{\Lambda}_0)=\Lambda_0+C'\delta$, for
some $C'\in\mathbb{Q}$.  \smallskip

Let $\mathbb{B}(\Lambda_{0})$ (resp.,
$\mathbb{B}\bigl(P_{\omega}^{\ast}(\check{\Lambda}_{0})\bigr)$) be
the set of all L-S paths of class $\Lambda_{0}$ (resp., of class
$P_{\omega}^{\ast}(\check{\Lambda}_{0})$). Let $\pi_{\Lambda_{0}}$
(resp., $\pi_{P_{\omega}^{\ast}(\check{\Lambda}_{0})}$) be the
straight path joining $0$ and $\Lambda_{0}$ (resp., $0$ and
$P_{\omega}^{\ast}(\check{\Lambda}_{0})$). For each integer $0\leq
i\leq e-1$, let $\widetilde{E}_i, \widetilde{F}_i$ denote the
raising root operator and the lowering root operator with respect
to the simple root $\alpha_i$.

\begin{lem} The map which sends $\pi_{P_{\omega}^{\ast}(\check{\Lambda}_0)}$
to $\pi_{\Lambda_{0}}$ extends to a bijection $\beta$ from
$\mathbb{B}\bigl(P_{\omega}^{\ast}(\check{\Lambda}_0)\bigr)$ onto
$\mathbb{B}(\Lambda_{0})$ such that $$
\beta\bigl(\widetilde{F}_{i_1}\cdots
\widetilde{F}_{i_s}\pi_{P_{\omega}^{\ast}(\check{\Lambda}_0)}\bigr)=\widetilde{F}_{i_1}\cdots
\widetilde{F}_{i_s}\pi_{\Lambda_{0}},
$$
for any $i_1,\cdots,i_s\in\{0,1,\cdots,e-1\}$.
\end{lem}

\noindent {Proof:} \,This follows from the fact that
$P_{\omega}^{\ast}(\check{\Lambda}_0)-\Lambda_{0}\in\mathbb{Q}\delta$
and the definitions of
$\mathbb{B}\bigl(P_{\omega}^{\ast}(\check{\Lambda}_0)\bigr)$ and
$\mathbb{B}(\Lambda_{0})$ (see \cite{Li1}).\medskip

Henceforth we shall identify
$\mathbb{B}\bigl(P_{\omega}^{\ast}(\check{\Lambda}_{0})\bigr)$
with $\mathbb{B}(\Lambda_{0})$. The action of $\omega^{\ast}$ on
$\mathfrak{h}^{\ast}$ naturally extends to the set
$\mathbb{B}\bigl(P_{\omega}^{\ast}(\check{\Lambda}_{0})\bigr)$
(and hence to the set $\mathbb{B}\bigl(\Lambda_{0}\bigr)$). By
\cite[(3.1.1)]{NS2}, if
$\widetilde{F}_{i_1}\widetilde{F}_{i_2}\cdots \widetilde{F}_{i_s}
\pi_{{\Lambda_{0}}}\in {\mathbb{B}}({\Lambda_{0}})$, then
\addtocounter{equation}{3}
\begin{equation}\label{equa66}
\omega^{\ast}\bigl(\widetilde{F}_{i_1}\widetilde{F}_{i_2}\cdots
\widetilde{F}_{i_s}
\pi_{{\Lambda_{0}}}\bigr)=\widetilde{F}_{\omega(i_{1})}\widetilde{F}_{\omega(i_{2})}\cdots
\widetilde{F}_{\omega(i_{s})}\pi_{{\Lambda_{0}}}.
\end{equation}

We denote by $\mathbb{B}^{\circ}\bigl(\Lambda_{0}\bigr)$ the set
of all L-S paths of class $\Lambda_{0}$ that are fixed by
$\omega^{\ast}$. For $\check{\mathfrak{g}}$, for each integer
$0\leq i\leq l$, we denote by $\widetilde{e}_i, \widetilde{f}_i$
the raising root operator and the lowering root operator with
respect to the simple root $\alpha_i$. Let
$\pi_{\check{\Lambda}_{0}}$ be the straight path joining $0$ and
$\check{\Lambda}_{0}$. By \cite[(4.2)]{NS1}, the linear map
$P_{\omega}^{\ast}$ naturally extends to a map from
$\check{\mathbb{B}}(\check{\Lambda}_{0})$ to
$\mathbb{B}^{\circ}\bigl(\Lambda_{0}\bigr)$ such that if
$\widetilde{f}_{i_1}\widetilde{f}_{i_2}\cdots \widetilde{f}_{i_s}
\pi_{\check{\Lambda}_{0}}\in
\check{\mathbb{B}}(\check{\Lambda}_{0})$, then (in the above two
cases)
$$\begin{aligned}
&P_{\omega}^{\ast}\bigl(\widetilde{f}_{i_1}\widetilde{f}_{i_2}\cdots
\widetilde{f}_{i_s}
\pi_{\check{\Lambda}_{0}}\bigr)=\omega\bigl(\widetilde{F}_{i_1}\bigr)
\omega\bigl(\widetilde{F}_{i_2}\bigr)\cdots
\omega\bigl(\widetilde{F}_{i_s}\bigr)\pi_{\Lambda_{0}},\end{aligned}
$$
where $$ \omega\bigl(\widetilde{F}_{i_t}\bigr):=\begin{cases}
\widetilde{F}_{i_t}\widetilde{F}_{\omega(i_{t})}, &\text{if
$c_{i_t,i_t}=2$ and $N_{i_t}=2$,}\\
\widetilde{F}_{i_{t}}, &\text{if $c_{i_t,i_t}=2$ and
$N_{i_t}=1$,}\\
\widetilde{F}_{\omega(i_t)}\widetilde{F}_{i_t}^2\widetilde{F}_{\omega(i_t)},
&\text{if $c_{i_t,i_t}=1$.}
\end{cases}
$$
Note that the case $c_{i_t,i_t}=1$ only happens when $e=2\ell+1$
and $i_t=\ell$.

\addtocounter{lem}{1}
\begin{lem} \label{lm35} {\rm(\cite[(4.2),(4.3)]{NS1})}
$\mathbb{B}^{\circ}\bigl(\Lambda_{0}\bigr)=P_{\omega}^{\ast}\bigl(
\check{\mathbb{B}}(\check{\Lambda}_{0})\bigr)$.
\end{lem}

Note that both $\check{\mathbb{B}}(\check{\Lambda}_{0})$ and
$\mathbb{B}\bigl(\Lambda_{0}\bigr)$ have a canonical crystal
structure with the raising and lowering root operators playing the
role of Kashiwara operators. They are isomorphic to the crystals
associated to the integrable highest weight modules
$\check{L}(\check\Lambda_{0})$ of highest weight
$\check\Lambda_{0}$
 over $U'_v(\check{\mathfrak{g}})$ and the integrable highest weight modules
 $L\bigl(\Lambda_{0}\bigr)$
of highest weight $\Lambda_{0}$ over $U'_v(\mathfrak{g})$
 respectively. Henceforth, we identify
them without further comments. Let $v_{\check{\Lambda}_{0}}$
(resp., $v_{\Lambda_{0}}$) denotes the unique highest weight
vector of highest weight $\check{\Lambda}_{0}$ (resp., of highest
weight $\Lambda_{0}$) in $\check{\mathbb{B}}(\check{\Lambda}_{0})$
(resp., in $\mathbb{B}(\Lambda_{0})$). Therefore, by
(\ref{equa66}) and Lemma \ref{lm35}, we get that

\addtocounter{cor}{5}
\begin{cor} \label{cor36}With the above assumptions and notations, there is an
injection $\eta$ from the set
$\check{\mathbb{B}}(\check{\Lambda}_{0})$ of crystal bases to the
set $\mathbb{B}(\Lambda_{0})$ of crystal bases such that
$$\begin{aligned}
&\eta\bigl(\widetilde{f}_{i_1}\widetilde{f}_{i_2}\cdots
\widetilde{f}_{i_s} v_{\check{\Lambda}_{0}}\bigr)\equiv
\omega\bigl(\widetilde{F}_{i_1}\bigr)\omega\bigl(\widetilde{F}_{i_2}\bigr)\cdots
\omega\bigl(\widetilde{F}_{i_s}\bigr)v_{\Lambda_{0}}\pmod{{v
L(\Lambda_{0})_{A}}},\end{aligned}
$$
where $i_1,\cdots,i_s$ are integers in $\{0,1,2,\cdots,\ell\}$,
and $A$ denotes the ring of rational functions in $\mathbb{Q}(v)$
which do not have a pole at $0$. Moreover, the image of $\eta$
consists of all crystal basis element $\widetilde{F}_{i_1}\cdots
\widetilde{F}_{i_s}v_{\Lambda_{0}}+v L(\Lambda_{0})_{A}$
satisfying $\widetilde{F}_{i_1}\cdots
\widetilde{F}_{i_s}v_{\Lambda_{0}}\equiv
\widetilde{F}_{\omega(i_{1})}\cdots
\widetilde{F}_{\omega(i_{s})}v_{\Lambda_{0}} \pmod{{v
L(\Lambda_{0})_{A}}}. $
\end{cor}

Let $\mathcal{K}:=\sqcup_{n\geq 0}\mathcal{K}_n$. We translate the
language of crystal bases into the language of partitions, we get
the following combinatorial result.

\addtocounter{thm}{6}
\begin{thm} \label{thm37} With the above notations, there is a
bijection $\eta$ from the set
$\check{\mathbb{B}}(\check{\Lambda}_{0})$ of crystal bases onto
the set $\bigl\{\lam\in\mathcal{K}\bigm|\MM(\lam)=\lam\bigr\}$,
such that if
$$
v_{\check{\Lambda}_{0}}\overset{r_1}{\twoheadrightarrow}\cdot
\overset{r_2}{\twoheadrightarrow}\cdot \cdots\cdots
\overset{r_s}{\twoheadrightarrow} \widetilde{f}_{r_s}\cdots
\widetilde{f}_{r_1}v_{\check{\Lambda}_{0}} $$ is a path from
$v_{\check{\Lambda}_{0}}$ to $\widetilde{f}_{r_s}\cdots
\widetilde{f}_{r_1}v_{\check{\Lambda}_{0}}$ in the crystal graph
of $L(\check{\lambda}_0)$, then the sequence
$$
\emptyset\underbrace{\overset{r_1}{\twoheadrightarrow}\cdot}_{\text{$\omega$
acts}}\,\,\underbrace{\overset{r_2}{\twoheadrightarrow}\cdot}_{\text{$\omega$
acts}}\cdot
\cdots\cdots\underbrace{\overset{r_s}{\twoheadrightarrow}\lam}_{\text{$\omega$
acts}}:=\eta\Bigl(\widetilde{f}_{r_s}\cdots
\widetilde{f}_{r_1}v_{\check{\Lambda}_{0}}\Bigr),
$$
where $$
\underbrace{\overset{r_t}{\twoheadrightarrow}\cdot}_{\text{$\omega$
acts}}:=\begin{cases} \overset{r_t}{\twoheadrightarrow}\cdot
\overset{e-r_t}{\twoheadrightarrow}, &\text{if $c_{r_t,r_t}=2$ and
$N_{r_t}=2$,}\\
\overset{r_t}{\twoheadrightarrow}, &\text{if $c_{r_t,r_t}=2$ and
$N_{r_t}=1$,}\\
\overset{\ell+1}{\twoheadrightarrow}\cdot\overset{\ell}{\twoheadrightarrow}\cdot
\overset{\ell}{\twoheadrightarrow}\cdot\overset{\ell+1}{\twoheadrightarrow}\cdot,
&\text{if $e=2\ell+1$ and $r_t=\ell$,}
\end{cases}
$$
defines a path in Kleshchev's $e$-good lattice which connects
$\emptyset$ and $e$-regular partition $\lam$ satisfying
$\MM(\lam)=\lam$.
\end{thm}

\noindent {Proof:} \,This follows from Lemma \ref{lm24}, Lemma
\ref{lm35} and Corollary \ref{cor36}.\medskip

For each partition $\lam$ of $n$, and each integer $0\leq i\leq
e-1$, we define
$$\begin{aligned}
\Sigma_i(\lam):&=\bigl\{\gamma\in[\lam]\bigm|\res(\gamma)=\overline{i}\bigr\},\\
N_i(\lam):&=\#\Sigma_i(\lam).\end{aligned} $$  Theorem 3.7 also
implies that if $\widetilde{f}_{r_1}\cdots
\widetilde{f}_{r_s}v_{\check{\Lambda}_{0}}\in\check{\mathbb{B}}(\check{\Lambda}_{0})$,
$\lam:=\eta\Bigl(\widetilde{f}_{r_1}\cdots
\widetilde{f}_{r_s}v_{\check{\Lambda}_{0}}\Bigr)$, then
\addtocounter{equation}{3}\begin{equation}\label{equa38}
N_i(\lam)=\begin{cases}\#\bigl\{1\leq t\leq s\bigm|r_t={i}\bigr\},
&\text{if $i\in\{0,1,2,\cdots,\ell-1\}$,}\\
\#\bigl\{1\leq t\leq s\bigm|r_t={e-i}\bigr\}, &\text{if
$i\in\{\ell+2,\ell+3,\cdots,e-1\}$,}\\
\#\bigl\{1\leq t\leq s\bigm|r_t={\ell-1}\bigr\}, &\text{if
$e=2\ell$ and $i=\ell+1$,}\\
\#\bigl\{1\leq t\leq s\bigm|r_t={\ell}\bigr\}, &\text{if
$e=2\ell$ and $i=\ell$,}\\
2\#\bigl\{1\leq t\leq s\bigm|r_t={\ell}\bigr\}, &\text{if
$e=2\ell+1$ and $i\in\{\ell,\ell+1\}$.}
\end{cases}
\end{equation}

\addtocounter{cor}{2}
\begin{cor} \label{cor39} Let $\lam\in\mathcal{K}_n$. Suppose that
$\MM(\lam)=\lam$.

1) If $e=2\ell+1$, then $N_{\ell}(\lam)=N_{\ell+1}(\lam)$.
Furthermore, $N_{\ell}(\lam)$ and $n-N_0(\lam)$ are both even
integers.

2) If $e=2\ell$, then $n-N_0(\lam)-N_{\ell}(\lam)$ is an even
integer.
\end{cor}

For each pair of integers $m, m'$ with $0\leq m+m'\leq n$, we
define
$$\begin{aligned}
\Sigma(n,m,m'):&=\bigl\{\lam\in\mathcal{K}_n\bigm|\MM(\lam)=\lam,
N_0(\lam)=m, N_{\ell}(\lam)=m'\bigr\},\\
N(n,m,m'):&=\#\Sigma(n,m,m'). \end{aligned}$$ Note that when
$e=2\ell+1$, by Corollary \ref{cor39}, $N(n,m,m')=0$ unless
$m+2m'\leq n$.

Recall the principle graduation introduced in \cite[\S1.5,
\S10.10]{Kac}. That is, the weight
$\Lambda_0-\sum_{i=0}^{e-1}k_i\alpha_i$ (where $k_i\in\Z$ for each
$i$) is assigned to degree $\sum_{i=0}^{e-1}k_i$. Let $\cch_t
L(\Lambda_0):=\sum_{n\geq 0}\dim L(\Lambda_0)_n t^n$ be the
principle specialized character\footnote{This is called
$q$-dimension in the book of Kac, see \cite[\S10.10]{Kac}.} of
$L(\Lambda_0)$, where
$L(\Lambda_0)_n=\oplus_{\deg\mu=n}L(\Lambda_0)_{\mu}$. Similarly,
let $L(\check{\Lambda}_0)$ denote the integrable highest weight
module of highest weight $\check{\Lambda}_0$ over
$\check{\mathfrak{g}}$. We use $\cch_t
L(\check{\Lambda}_0):=\sum_{n\geq 0}\dim L(\check{\Lambda}_0)_n
t^n$ to denote the principle specialized character of
$L(\check{\Lambda}_0)$. Now applying Lemma \ref{lm24}, Lemma
\ref{lm35} and Theorem \ref{thm37}, we get that
\addtocounter{equation}{1}\begin{equation}\label{equa310}\dim
L(\check{\Lambda}_0)_{n}=\sum_{0\leq m+m'\leq n}N(2n-m+2m',m,2m').
\end{equation}
if $e=2\ell+1$; while
\begin{equation}\label{equa311}\dim
L(\check{\Lambda}_0)_{n}=\sum_{0\leq m+m'\leq n}N(2n-m-m',m,m').
\end{equation}
if $e=2\ell$.\medskip

Suppose that $e=2\ell+1$. That is, we are in the odd case. In this
case, $\check{\mathfrak{g}}$ is the twisted affine Lie algebra of
type $A_{2\ell}^{(2)}$. By \cite[(14.5.4)]{Kac}, the principle
specialized character of $L(\check{\Lambda}_0)$ is given by
\begin{equation}\label{equa312}
\cch_t L(\check{\Lambda}_0)=\prod_{\substack{\text{$i\geq 1$, $i$ odd}\\
\text{$i\not\equiv 0\!\!\!\pmod{e}$ }}}\frac{1}{1-t^i}.
\end{equation}
Hence by (\ref{equa310}) and (\ref{equa312}), we get that

\addtocounter{thm}{5} \begin{thm} \label{thm313} With the above
notations, we have that
$$
\prod_{\substack{\text{$i\geq 1$, $i$ odd}\\
\text{$i\not\equiv 0\!\!\!\pmod{e}$ }}}\frac{1}{1-t^i}=
\sum_{n\geq 0}\biggl(\sum_{0\leq m+m'\leq
n}N(2n-m+2m',m,2m')\biggr)t^n.$$
\end{thm}

In \cite{K}, Kang has given a combinatorial realization of
$\check{\mathbb{B}}(\check{\Lambda}_{0})$ in terms of reduced
proper Young walls, which are inductively defined. In our
$A_{2\ell}^{(2)}$ case, a direct explicit characterization can be
given in terms of restricted $e$-strict partitions as follows, see
\cite{LT2},\cite{BK1}.

Recall that (\cite{BK1},\cite{BK2}) a partition $\lam$ is called
$e$-strict if $\lam_i=\lam_{i+1}\,\Rightarrow\,e|\lam_i$ for each
$i=1,2,\cdots $. An $e$-strict partition $\lam$ is called
restricted if in addition $$
\begin{cases}\lam_i-\lam_{i+1}\leq e, &\text{if $e\nmid\lam_i$,}\\
\lam_i-\lam_{i+1}<e, &\text{if $e|\lam_i$.}
\end{cases}
\quad\text{for each $i=1,2,\cdots.$}
$$
Let $DPR_e(n)$ denote the set of all restricted $e$-strict
partitions of $n$. Let $DPR_e:=\sqcup_{n\geq 0}DPR_e(n)$.

It turns out that there is a natural $1$-$1$ correspondence
between $\check{\mathbb{B}}(\check{\Lambda}_{0})$ and $DPR_e$.
Furthermore, the crystal structure
$\check{\mathbb{B}}(\check{\Lambda}_{0})$ can be concretely
realized via some combinatorics of $DPR_e$, which we now describe.

We recall some notions. Elements of $(r,s)\in\Z_{>0}\times
\Z_{>0}$ are called nodes. Let $\lam$ be an $e$-strict partition.
We label the nodes of $\lam$ with {\it residues}, which are the
elements of $\Z/(\ell+1)\Z$. The residue of the node $A$ is
denoted $\res{A}$. The labelling depends only on the column and
following the repeating pattern $$
\overline{0},\overline{1},\cdots,
\overline{\ell-1},\overline{\ell},\overline{\ell-1},\cdots,\overline{1},\overline{0},
$$
starting from the first column and going to the right. For
example, let $e=5$, $\ell=2$, let $\lam=(10,10,6,1)$ be a
restricted $5$-strict partition of $27$. Its residues are as
follows: $$
\begin{matrix}
\overline{0} & \overline{1} &\overline{2} &\overline{1}
&\overline{0} &\overline{0} &\overline{1} &\overline{2} &
\overline{1} & \overline{0}\\
\overline{0} & \overline{1} &\overline{2} &\overline{1}
&\overline{0} &\overline{0} &\overline{1} &\overline{2} &
\overline{1} & \overline{0}\\
\overline{0} & \overline{1} &\overline{2} &\overline{1}
&\overline{0} &\overline{0} & & &
& \\
\overline{0} & & & & &  & & &\end{matrix}
$$
A node $A=(r,s)\in[\lam]$ is called removable (for $\lam$) if
either

R1) $\lam_{A}:=\lam-\{A\}$ is again an $e$-strict partition; or

R2) the node $B=(r,s+1)$ immediately to the right of $A$ belongs
to $\lam$, $\res(A)=\res(B)$, and both $\lam_B$ and
$\lam_{AB}:=\lam-\{A,B\}$ are $e$-strict partitions.\smallskip

Similarly, a node $B=(r,s)\notin[\lam]$ is called addable (for
$\lam$) if either

A1) $\lam^{B}:=\lam\cup\{B\}$ is again an $e$-strict partition; or

A2) the node $A=(r,s-1)$ immediately to the left of $B$ does not
belong to $\lam$, $\res(A)=\res(B)$, and both
$\lam^A:=\lam\cup\{A\}$ and $\lam^{AB}:=\lam\cup\{A,B\}$ are
$e$-strict partitions.

Note that R2) and A2) above are only possible for nodes with
residue $\overline{0}$. Now fix a residue $x$ and consider the
sequence of removable and addable $x$-nodes obtained by reading
the boundary of $\lam$ from the bottom left to top right. We use
``A'' to denote an addable $x$-node and use ``R'' to denote a
removable $x$-node, then we get a sequence of letters A,R. Given
such a sequence, we remove all occurrences of the string ``AR''
and keep on doing this until no such string ``AR'' is left. The
``R''s that still remain are the {\it normal} $x$-nodes of $\lam$
and the rightmost of these is the {\it good} $x$-node, the ``A''s
that still remain are the {\it conormal} $x$-nodes of $\lam$ and
the leftmost of these is the {\it cogood} $x$-node. Note
that\footnote{This is because any removable node $\gamma$ of type
R2) has an adjacent neighborhood $\gamma'$ in his right, which is
another removable node with the same residue. If $\gamma$ could
survive after deleting all the string ``AR'', then $\gamma'$ must
also survive. In that case, $\gamma'$ is a normal node higher than
$\gamma$. So $\gamma$ can not be a good node. For cogood node of
type A2), the reason is similar.} good $x$-node is necessarily of
type R1), and cogood $x$-node is necessarily of type A1). We
define
$$\begin{aligned} \eps_i(\lam)&=\#\bigl\{\text{$i$-normal nodes in
$\lam$}\bigr\},\\
\varphi_i(\lam)&=\#\bigl\{\text{$i$-conormal nodes in
$\lam$}\bigr\}\end{aligned}
$$
and we set $$\begin{aligned}
\widetilde{e}_{i}(\lam)&=\begin{cases} \lam_{A},
&\text{if $\eps_i(\lam)>0$ and $A$ is the (unique) good $i$-node,}\\
0, &\text{if $\eps_i(\lam)=0$.}
\end{cases}\\
\widetilde{f}_{i}(\lam)&=\begin{cases} \lam^{B},
&\text{if $\varphi_i(\lam)>0$ and $B$ is the (unique) cogood $i$-node,}\\
0, &\text{if $\varphi_i(\lam)=0$.}
\end{cases}
\end{aligned}$$ Then, we get an infinite colored oriented graph, whose vertices are
$e$-strict partitions and whose arrows are given by $$
\text{$\lam\overset{i}{\twoheadrightarrow}\mu$\,\,\,$\Longleftrightarrow$\,
$\mu=\widetilde{f}_{i}(\lam)$\,\,$\Longleftrightarrow$\,\,$\lam=\widetilde{e}_{i}(\mu)$}.
$$

The sublattice spanned by all restricted $e$-strict partitions
equipped with the functions $\varepsilon_i, \varphi_i$ and the
operators $\widetilde{e}_i, \widetilde{f}_i$, can be turned into a
colored oriented graph which we denote by ${\mathfrak{RP}}_{e}$.

\addtocounter{lem}{8}
\begin{lem} \label{lem314}{\rm (\cite{K})} With the above notations, the
graph $\mathfrak{RP}_e$ can be identified with the crystal graph
$\check{\mathbb{B}}(\check{\Lambda}_{0})$ associated to the
integrable highest weight $\check{\mathfrak{g}}$-module of highest
weight $\check{\Lambda}_0$.
\end{lem}

Applying Theorem \ref{thm37} and Lemma \ref{lem314}, we get that
\addtocounter{thm}{1}
\begin{thm} \label{thm315} With the above notations, there is a
bijection $\eta$ from the set $DPR_e$ of restricted $e$-strict
partitions onto the set
$\bigl\{\lam\in\mathcal{K}\bigm|\MM(\lam)=\lam\bigr\}$, such that
if
$$
\emptyset\overset{r_1}{\twoheadrightarrow}\cdot
\overset{r_2}{\twoheadrightarrow}\cdot \cdots\cdots
\overset{r_s}{\twoheadrightarrow} \check{\lam} $$ is a path from
$\emptyset$ to $\check{\lam}$ in the subgraph $\mathfrak{RP}_e$,
then the sequence then the sequence
$$
\emptyset\underbrace{\overset{r_1}{\twoheadrightarrow}\cdot
\overset{2\ell+1-r_1}{\twoheadrightarrow}\cdot
}_{\text{$\widetilde{N_{r_1}}$
terms}}\underbrace{\overset{r_2}{\twoheadrightarrow}\cdot
\overset{2\ell+1-r_2}{\twoheadrightarrow}\cdot}_{\text{$\widetilde{N_{r_2}}$
terms}}
\cdots\cdots\underbrace{\overset{r_s}{\twoheadrightarrow}\cdot
\overset{2\ell+1-r_s}{\twoheadrightarrow}\lam}_{\text{$\widetilde{N_{r_s}}$
terms}}:=\eta\bigl(\check{\lam}\bigr),
$$
where $$ \underbrace{\overset{r_t}{\twoheadrightarrow}\cdot
\overset{2\ell+1-r_t}{\twoheadrightarrow}\cdot}_{\text{$\widetilde{N_{r_t}}$
terms}}:=\begin{cases}
\overset{r_t}{\twoheadrightarrow}\cdot\overset{2\ell+1-r_t}{\twoheadrightarrow}\cdot
&\text{if
$r_t\in\{1,2,\cdots,\ell-1\}$,}\\
\overset{0}{\twoheadrightarrow}\cdot, &\text{if
$r_t=0$,}\\
\overset{\ell+1}{\twoheadrightarrow}\cdot
\overset{\ell}{\twoheadrightarrow}\cdot
\overset{\ell}{\twoheadrightarrow}\cdot
\overset{\ell+1}{\twoheadrightarrow}\cdot, &\text{if $r_t=\ell$,}
\end{cases} $$
defines a path in Kleshchev's $(2\ell+1)$-good lattice which
connects $\emptyset$ and $(2\ell+1)$-regular partition $\lam$
satisfying $\MM(\lam)=\lam$.
\end{thm}

\noindent {\bf Remark 3.16}\,\,\,In \cite{BK1}, \cite{BK2},
Brundan and Kleshchev investigated the modular representations of
Hecke-Clifford superalgebras at defining parameter a primitive
$(2\ell+1)$-th root of unity as well as of affine Sergeev
superalgebras over a field of characteristic $2\ell+1$. Their main
result states that the modular socle branching rules of these
superalgebras provide a realization of the crystal of the twisted
affine Lie algebra of type $A_{2\ell}^{(2)}$. This applies, in
particular, to the modular socle branching rules of the spin
symmetric group $\widehat{\BS}_n$, which is the double cover of
the symmetric group $\BS_n$. It would be interesting to know if
there is any connection between their results and ours, at least
in the special case where $q=1$ and $2\ell+1$ being a prime
number.
\medskip

Let $P_n$ be the set of all partitions of $n$. Let
$P:=\sqcup_{n\geq 0}P_n$. Recall that when $\HH_k(\BS_n)$ is
semisimple, then $\mathcal{K}_n=P_n$ and $\MM$ degenerates to the
map $\lam\mapsto\lam^t$ for any $\lam\in P_n$. Let $DP_n$ be the
set of all partitions into distinct parts (i.e., the set all
$0$-strict partitions). Let $DP:=\sqcup_{n\geq 0}DP_n$. Let $SP$
be the set of all symmetric partitions, i.e., $SP:=\bigl\{\lam\in
P \bigm|\lam=\lam^t\bigr\}$. We shall now establish a bijection
between the set $DP$ and the set $SP$. Note that in the special
case where $q=1$ and $2\ell+1$ is a prime number, the set $DP_n$
parameterizes the ordinary irreducible supermodules of the spin
symmetric group $\widehat{\BS}_n$, while the set
$SP_n:=\bigl\{\lam\in P_n \bigm|\lam=\lam^t\bigr\}$ parameterizes
those ordinary irreducible modules of the symmetric group which
splits on restriction to the alternating group $A_n$.

For each partition $\lam=(\lam_1,\lam_2,\cdots,\lam_s)\in DP$ with
$\ell(\lam)=s$, let
$\lam^t=(\lam^t_1,\lam^t_2,\cdots,\lam^t_{\lam_1})$ be the
conjugate of $\lam$, we define
$$
\widetilde{\eta}\bigl(\lam\bigr)=(\lam_1,\lam_2+1,\lam_3+2,\cdots,\lam_s+s-1,\lam^{t}_{s+1},
\lam^{t}_{s+2},\cdots,\lam^t_{\lam_1}).
$$

\addtocounter{thm}{1}
\begin{thm} \label{thm317} With the above notations, the map $\widetilde{\eta}$
defines a bijection from the set $DP$ onto the set $SP$.
\end{thm}

\noindent {Proof:} \,Let $\lam\in DP$. By definition,
$\lam_1>\lam_2>\cdots>\lam_s$, it follows that $$
\lam_1\geq\lam_2+1\geq\lam_3+2\geq\cdots\geq\lam_s+s-1\geq\lam^{t}_{s+1}\geq
\lam^{t}_{s+2}\geq\cdots\geq\lam^t_{\lam_1}.
$$
That is, $\widetilde{\eta}\bigl(\lam\bigr)\in P$. We claim that
$\Bigl(\widetilde{\eta}\bigl(\lam\bigr)\Bigr)^t=\widetilde{\eta}\bigl(\lam\bigr)$.

We use induction on $\ell(\lam)$. Suppose that
$\Bigl(\widetilde{\eta}\bigl(\nu\bigr)\Bigr)^t=\widetilde{\eta}\bigl(\nu\bigr)$
for any partition $\nu$ satisfying $\ell(\nu)<\ell(\lam)$. We
write
$\mu=(\mu_1,\cdots,\mu_{\lam_1})=\widetilde{\eta}\bigl(\lam\bigr)$.
Then $$ \mu_i=\begin{cases} \lam_i+i-1, &\text{for $1\leq i\leq
s$,}\\ \lam^t_i, &\text{for $s+1\leq i\leq\lam_1$.}
\end{cases}\quad\text{for $i=1,2,\cdots,\lam_1$.}
$$
By definition, $\mu_i^t=\#\{1\leq j\leq\lam_1|\mu_j\geq i\}$. It
is clear that $\mu^t_1=\lam_1=\mu_1$. We remove away the first row
as well as the first column of $\mu$. Then we get a partition
$\widehat{\mu}$. It is easy to see that $$
\widehat{\mu}=(\lam_2,\lam_3+1,\cdots,\lam_s+s-2,\lam^{t}_{s+1}-1,
\lam^{t}_{s+2}-1,\cdots,\lam^t_{\lam_2}-1)=\widetilde{\eta}\bigl(\widehat{\lam}\bigr),
$$
where $\widehat{\lam}:=(\lam_2,\lam_3,\cdots,\lam_s)$.

Note that $\ell(\widehat{\lam})<\ell(\lam)$. By induction
hypothesis, we know that
$\bigl(\widehat{\mu}\bigr)^t=\widehat{\mu}$. It follows that
$\mu^t=\mu$ as well. This proves our claim.\smallskip

Second, we claim that the map $\widetilde{\eta}$ is injective. In
fact, suppose that $$\begin{aligned}
\widetilde{\eta}\bigl(\lam\bigr)&=(\lam_1,\lam_2+1,\lam_3+2,\cdots,\lam_s+s-1,\lam^{t}_{s+1},
\lam^{t}_{s+2},\cdots,\lam^t_{\lam_1})\\
&=(\mu_1,\mu_2+1,\mu_3+2,\cdots,\mu_{s'}+s'-1,\mu^{t}_{s'+1},
\mu^{t}_{s'+2},\cdots,\mu^t_{\mu_1})=
\widetilde{\eta}\bigl(\mu\bigr),\end{aligned}$$ where $\lam,\mu\in
DP$, $\ell(\lam)=s, \ell(\mu)=s', s\leq s'$. Then $$ \lam_1
=\ell\Bigl(\widetilde{\eta}\bigl(\lam\bigr)\Bigr)=\ell\Bigl(\widetilde{\eta}\bigl(\mu\bigr)\Bigr)
=\mu_1.
$$
It follows that $\lam_i=\mu_i$ for $i=1,2,\cdots,s$. If $s<s'$,
then $\lam_{s+1}^t=\mu_{s+1}+s\geq s+1$, which is impossible.
Therefore $s=s'$, and hence $\lam=\mu$. This proves the
injectivity of $\widetilde{\eta}$.

It remains to show that $\widetilde{\eta}$ is surjective. Let
$\mu\in P$ such that $\mu^t=\mu$. Let $A=(r,s)$ be the unique node
on the boundary of $[\lam]$ which sits on the main diagonal of
$[\lam]$. We define
$$ \lam:=(\mu_1,\mu_2-1,\mu_3-2,\cdots,\mu_{r}-r+1).
$$
Then one sees easily that $\lam\in DP$ and
$\widetilde{\eta}(\lam)=\mu$. This proves that $\widetilde{\eta}$
is surjective, hence completes the proof of the whole theorem.
\medskip

\noindent {\bf Remark 3.18}\,\,\,We remark that if one consider
the special case where $q=1$ and $2\ell+1$ is a prime number, it
would be interesting to know if the reduced decomposition matrices
(in the sense of \cite[(6.2)]{LT2}) of the spin symmetric groups
are embedded as submatrices into the decomposition matrices of the
alternating groups in odd characteristic $e$ via our bijections
$\eta$ and $\widetilde{\eta}$.
\bigskip

Now we suppose that $e=2\ell$. That is, we are in the even case.
In this case, $\check{\mathfrak{g}}$ is the twisted affine Lie
algebra of type $D_{\ell+1}^{(2)}$. By \cite[(14.5.4)]{Kac}, the
principle specialized character of $L(\check{\Lambda}_0)$ is given
by \addtocounter{equation}{6}
\begin{equation}\label{equa319}
\cch_t L(\check{\Lambda}_0)=\prod_{\text{$i\geq 1$, $i$
odd}}\frac{1}{1-t^i}.
\end{equation}
Hence by (\ref{equa311}) and (\ref{equa319}), we get that

\addtocounter{thm}{2} \begin{thm} \label{thm320} With the above
notations, we have that
$$
\prod_{\text{$i\geq 1$, $i$ odd}}\frac{1}{1-t^i}= \sum_{n\geq
0}\biggl(\sum_{0\leq m+m'\leq n}N(2n-m-m',m,m')\biggr)t^n.$$
\end{thm}

We propose the following definition.

\addtocounter{dfn}{20}
\begin{dfn} \label{dfn321} Let $f\in\mathbb{N}$ with $f>1$. An $f$-strict partition $\lam$
is called double restricted if
$$\begin{cases}\lam_i-\lam_{i+1}\leq 2f, &\text{if $f\nmid\lam_i$,}\\
\lam_i-\lam_{i+1}<2f, &\text{if $f|\lam_i$.}
\end{cases}
\quad\text{for each $i=1,2,\cdots.$}
$$
Here we make the convention that $\lam_i=0$ for any
$i>\ell(\lam)$.
\end{dfn}
Let $DDPR_f(n)$ denote the set of all double restricted $f$-strict
partitions of $n$. Let $DDPR_f:=\sqcup_{n\geq
0}DDPR_f(n)$.\smallskip

In \cite{K}, Kang has given a combinatorial realization of
$\check{\mathbb{B}}(\check{\Lambda}_{0})$ in terms of reduced
proper Young walls, which are inductively defined. In our
$D_{\ell+1}^{(2)}$ case, we shall give a direct explicit
characterization in terms of double restricted $(\ell+1)$-strict
partitions as follows.

As before, elements of $(r,s)\in\Z_{>0}\times \Z_{>0}$ are called
nodes. Let $\lam$ be an $(\ell+1)$-strict partition. We label the
nodes of $\lam$ with {\it residues}, which are the elements of
$\Z/(\ell+1)\Z$. The residue of the node $A$ is denoted $\res{A}$.
The labelling depends only on the column and following the
repeating pattern $$ \overline{0},\overline{1},\cdots,
\overline{\ell-1},\overline{\ell},\overline{\ell},\overline{\ell-1},\cdots,\overline{1},\overline{0},
$$
starting from the first column and going to the right. For
example, let $e=4$, $\ell=2$, let $\lam=(9,9,7,1)$ be a double
restricted $3$-strict partition of $26$. Its residues are as
follows: $$
\begin{matrix}
\overline{0} & \overline{1} & \overline{2} & \overline{2} &
\overline{1} & \overline{0} & \overline{0} & \overline{1} &
\overline{2} \\
\overline{0} & \overline{1} & \overline{2} & \overline{2} &
\overline{1} & \overline{0} & \overline{0} & \overline{1} &
\overline{2} \\
\overline{0} & \overline{1} & \overline{2} & \overline{2}
& \overline{1} & \overline{0} & \overline{0} & &  \\
\overline{0} & & & & & & & &
\end{matrix}
$$

Let $\lam$ be an $(\ell+1)$-strict partition. A node
$A=(r,s)\in[\lam]$ is called removable (for $\lam$) if either

R1) $\lam_{A}:=\lam-\{A\}$ is again an $(\ell+1)$-strict
partition; or

R2) the node $B=(r,s+1)$ immediately to the right of $A$ belongs
to $\lam$, $\res(A)=\res(B)$, and both $\lam_B$ and
$\lam_{AB}:=\lam-\{A,B\}$ are $(\ell+1)$-strict
partitions.\smallskip

Similarly, a node $B=(r,s)\notin[\lam]$ is called addable (for
$\lam$) if either

A1) $\lam^{B}:=\lam\cup\{B\}$ is again an $(\ell+1)$-strict
partition; or

A2) the node $A=(r,s-1)$ immediately to the left of $B$ does not
belong to $\lam$, $\res(A)=\res(B)$, and both
$\lam^A:=\lam\cup\{A\}$ and $\lam^{AB}:=\lam\cup\{A,B\}$ are
$(\ell+1)$-strict partitions. Now we can define the notions of
normal (resp., comormal) nodes, good (resp., cogood) nodes, the
functions $\varepsilon_i, \varphi_i$ and the operators
$\widetilde{e}_i, \widetilde{f}_i$ in the same way as in the case
where $e=2\ell+1$. Note that the definition of residue in the even
case is different with the odd case, and in the even case we deal
with $(\ell+1)$-strict partitions instead of $e$-strict
partitions.

\addtocounter{lem}{7}
\begin{lem} \label{lem322} Let $\lam$ be any given double restricted $(\ell+1)$-strict
partition. Then

1) there exists good (removable) node as well as cogood (addable)
node for $\lam$;

2) for any good (removable) node $A$ for $\lam$, $\lam-\{A\}$ is
again a double restricted $(\ell+1)$-strict partition. In
particular, there is a path (not necessary unique) from the empty
partition $\emptyset$ to $\lam$ in the lattice spanned by double
restricted $(\ell+1)$-strict partitions;

3) for any cogood (addable) node $A$ for $\lam$, $\lam\cup\{A\}$
is again a double restricted $(\ell+1)$-strict partition.
\end{lem}

\noindent {Proof:} \,We write $\lam=(\lam_1,\cdots,\lam_s)$, where
$\ell(\lam)=s$. Let $B=(s,\lam_s)$. Then, as $\lam$ is double
restricted, either $\lam_s=1$ or $\lam_s>1$ and
$\res(B)\neq\overline{0}$. In both cases, one sees easily that $B$
must be a normal $\res(B)$-node (as there are no addable
$\res(B)$-nodes below $B$). It follows that there must exist good
(removable) $\res(B)$-node for $\lam$. In a similar way, one can
show that $B'=(1,\lam_1+1)$ is a conormal $\res(B')$-node, which
implies that there must exist cogood (addable) $\res(B')$-node for
$\lam$. This proves 1).

Now let $A=(a,\lam_a)$ be a good (removable) node for $[\lam]$.
Then $A$ is necessarily of type R1). If $a=1$, then it is easy to
check that $\lam-\{A\}$ is again double restricted
$(\ell+1)$-strict. Suppose that $a>1$. We write $\res(A)=i$. We
claim that $\lam_{a-1}-\lam_a<2(\ell+1)$. In fact, If
$\lam_{a-1}-\lam_a=2(\ell+1)$, then either $\lam_a\not\equiv
0\!\pmod{\ell+1}$, or $\lam_a\equiv 0\!\pmod{\ell+1}$. In the
former case, one sees easily that $(a-1,\lam_{a-1})$ is a
removable $i$-node of type R1) next to (the right of) $A$ and there is no addable 
$i$-node sitting between them. Now as $A$ survives after deleting all the string ``AR'', the
node $(a-1,\lam_{a-1})$ must also survive after deleting all the
string ``AR''. In other words, it is in fact a normal $i$-node of
$\lam$ higher than $A$, which is impossible (since $A$ is the
unique good $i$-node of $\lam$); while in the latter case, it
would follows that $\lam_{a-1}\equiv 0\!\pmod{\ell+1}$, and hence
$\lam_{a-1}-\lam_a<2(\ell+1)$ because $\lam$ is double restricted
$(\ell+1)$-strict, which is again a contradiction. This proves our
claim. Now there are only five possibilities:
\smallskip

\noindent {\it Case
1.}\,\,$i\notin\{\overline{0},\overline{\ell}\}$.\smallskip

Then either $\lam_{a-1}\not\equiv 0\!\pmod{\ell+1}$ or
$\lam_{a-1}\equiv 0\!\pmod{\ell+1}$ and
$\lam_{a-1}-\lam_a<2\ell+1$. In both cases, one checks easily that
$\lam-\{A\}$ is again a double restricted $(\ell+1)$-strict.
\smallskip

\noindent {\it Case 2.}\,\,$i=\overline{\ell}$ and $\lam_a\equiv
0\!\pmod{\ell+1}$.\smallskip

Since $\lam_{a-1}-\lam_a<2(\ell+1)$, it follows that
$\lam_{a-1}-(\lam_a-1)\leq 2(\ell+1)$. Now $\lam_a\equiv
0\!\pmod{\ell+1}$ implies that either $\lam_{a-1}\not\equiv
0\!\pmod{\ell+1}$ or $\lam_{a-1}=\lam_a+\ell+1$. In both cases one
sees easily that $\lam$ is double restricted $(\ell+1)$-strict
must imply that $\lam-\{A\}$ is double restricted
$(\ell+1)$-strict too.
\smallskip

\noindent {\it Case 3.}\,\,$i=\overline{\ell}$ and $\lam_a\equiv
1\!\pmod{\ell+1}$.\smallskip

We know that $\lam_{a-1}-\lam_a<2(\ell+1)$. We claim that
$\lam_{a-1}-\lam_a<2\ell+1$. In fact, if $\lam_{a-1}-\lam_a=
2\ell+1$, then $(a-1,\lam_{a-1})$ must be another normal
$\overline{\ell}$-node higher than $A$, which is impossible. This
proves our claim. Therefore, $\lam_{a-1}-(\lam_a-1)\leq 2\ell+1$,
which implies that $\lam-\{A\}$ is still double restricted
$(\ell+1)$-strict.\smallskip

\noindent {\it Case 4.}\,\,$i=\overline{0}$ and $\lam_a\equiv
0\!\pmod{2(\ell+1)}$.\smallskip

In this case one proves that $\lam-\{A\}$ is double restricted
$(\ell+1)$-strict by using the same argument as in the proof of
Case 2.
\smallskip

\noindent {\it Case 5.}\,\,$i=\overline{0}$ and $\lam_a\equiv
1\!\pmod{2(\ell+1)}$.\smallskip

In this case one proves that $\lam-\{A\}$ is double restricted
$(\ell+1)$-strict by using the same argument as in the proof of
Case 3.

This completes the proof of 2). The proof of 3) is similar and is
left to the readers.
\medskip

Therefore, the lattice spanned by all double restricted
$(\ell+1)$-strict partitions equipped with the functions
$\varepsilon_i, \varphi_i$ and the operators $\widetilde{e}_i,
\widetilde{f}_i$, can be turned into a colored oriented graph
which we denote by $\widetilde{\mathfrak{RP}}_{\ell+1}$.

\begin{lem} The graph $\widetilde{\mathfrak{RP}}_{\ell+1}$
can be identified with the crystal graph
$\check{\mathbb{B}}(\check{\Lambda}_{0})$ associated to the
integrable highest weight $\check{\mathfrak{g}}$-module of highest
weight $\check{\Lambda}_0$.\end{lem}

\noindent {Proof:} \,This follows from Lemma \ref{lem322} and
Kang's combinatorial construction of the proper Young wall (see
\cite{K} and \cite{HK}). Note that our definition of removable and
addable node are in accordance with the definition given in
\cite[page 275, 278]{K}. To translate the language of proper Young
walls into the language of double restricted strict partitions,
one has to think the columns of the Young walls in \cite{K} as the
rows of our double restricted strict partitions. \smallskip

Applying Theorem \ref{thm37}, we get that

\addtocounter{thm}{3}
\begin{thm} \label{thm324} With the above notations, there is a
bijection $\eta$ from the set $DDPR_{\ell+1}$ of double restricted
$(\ell+1)$-strict partitions onto the set
$\bigl\{\lam\in\mathcal{K}\bigm|\MM(\lam)=\lam\bigr\}$, such that
if
$$
\emptyset\overset{r_1}{\twoheadrightarrow}\cdot
\overset{r_2}{\twoheadrightarrow}\cdot \cdots\cdots
\overset{r_s}{\twoheadrightarrow} \check{\lam} $$ is a path from
$\emptyset$ to $\check{\lam}$ in the graph
$\widetilde{\mathfrak{RP}}_{\ell+1}$, then the sequence
$$
\emptyset\underbrace{\overset{r_1}{\twoheadrightarrow}\cdot
\overset{2\ell-r_1}{\twoheadrightarrow}\cdot}_{\text{$N_{r_1}$
terms}}\underbrace{\overset{r_2}{\twoheadrightarrow}\cdot
\overset{2\ell-r_2}{\twoheadrightarrow}\cdot}_{\text{$N_{r_2}$
terms}} \cdots\cdots\underbrace{\cdot
\overset{r_s}{\twoheadrightarrow}\cdot
\overset{2\ell-r_s}{\twoheadrightarrow}\lam}_{\text{$N_{r_s}$
terms}}:=\eta\bigl(\check{\lam}\bigr),
$$
where $$ \underbrace{\overset{r_t}{\twoheadrightarrow}\cdot
\overset{2\ell-r_t}{\twoheadrightarrow}\cdot}_{\text{$N_{r_t}$
terms}}:=\begin{cases} \overset{r_t}{\twoheadrightarrow}\cdot
\overset{2\ell-r_t}{\twoheadrightarrow}\cdot, &\text{if
$r_t\in\{1,2,\cdots,\ell-1\}$,}\\
 \overset{r_t}{\twoheadrightarrow}\cdot, &\text{if
$r_t\in\{{0},{\ell}\}$,}\\
\end{cases}
$$ defines a path in Kleshchev's $(2\ell)$-good lattice which
connects $\emptyset$ and $(2\ell)$-regular partition $\lam$
satisfying $\MM(\lam)=\lam$.
\end{thm}

\noindent {\bf Remark 3.25}\,\,\,In \cite{LT2}, Leclerc-Thibon
conjectured that the decomposition matrices of Hecke-Clifford
superalgebras with parameter $q$ should related to the Fock space
representation of the twisted affine Lie algebra of type
$A_{2\ell}^{(2)}$ if $q$ is a primitive $(2\ell+1)$th root of
unity; or of type $D_{\ell+1}^{(2)}$ if $q$ is a primitive
$2\ell$-th root of unity. In \cite{BK1}, \cite{BK2}, Brundan and
Kleshchev show that the modular irreducible super-representations
of Hecke-Clifford superalgebras at defining parameter $q$ a
primitive $(2\ell+1)$-th root of unity as well as of affine
Sergeev superalgebras over a field of characteristic $2\ell+1$ are
parameterized by the set of restricted $(2\ell+1)$-strict
partitions, which partly verified the idea of \cite{LT2}. It would
be interesting to know if our notion of double restricted
$(\ell+1)$-strict partitions give a natural parameterization of
the modular irreducible super-representations of Hecke-Clifford
superalgebras when $q$ is a primitive $(2\ell)$-th root of unity.

\bigskip\bigskip
\bigskip \bigskip
\centerline{ACKNOWLEDGEMENT} \bigskip

\thanks{Research supported by the URF of Victoria University of
Wellington and the National Natural Science Foundation of China
(Project 10401005) and the Program NCET. The author wishes to
thank the School of Mathematics, Statistics and Computer Science
at Victoria University of Wellington for their hospitality during
his visit in 2005. He also appreciates the referee for several
helpful comments and for pointing out an error in the first
version of this paper.}

\end{document}